\documentclass[11pt]{amsart}
\usepackage[T1]{fontenc}
\usepackage{latexsym,amssymb,amsmath}
\input epsf

\def\vs{variété de Schubert\;}
\def\vss{variétés de Schubert\;}

\def\cR{{\mathcal R}}
\def\cS{{\mathcal S}}

\def\11{\mathbf 1}

\def\FF{\mathbb F}

\def\e{\varepsilon}

\def\dim{{\rm dim}\;}

\newcommand\dem{{\noindent {\em Démonstration}.}\hspace{2mm}}
\newtheorem{theo}{Théorème}

\newtheorem{lemm}[theo]{Lemme}
\newtheorem{prop}[theo]{Proposition}

\begin{document}

\title{Le lieu singulier des variétés de Schubert}
\author{L. Manivel}
\date{February 2001}

\medskip
\begin{abstract}
On donne une description explicite des composantes irréductibles du
lieu singulier d'une variété de Schubert dans une variété de drapeaux
classique. Cette description implique une conjecture de Lakshmibai et
Sandhya. 
\end{abstract}
\maketitle

%\medskip
\section{Enoncé du résultat}

On se place sur un corps algébriquement clos de caractéristique
nulle. La variété $\FF_n$ des drapeaux complets de sous-espaces 
d'un espace vectoriel de dimension $n$ est une variété projective
homogène de dimension $n(n-1)/2$. Elle admet une décomposition
cellulaire classique dont les adhérences, les variétés de Schubert,
sont indexées par les permutations de l'intervalle $[1,\ldots ,n]$.
Ces variétés dépendent du choix d'un drapeau de référence
$V_{\bullet}$, en fonction desquelles elles se définissent de la façon
suivante : si $w \in\cS_n$,
$$X_w=\{W_{\bullet}\in \FF_n,
\;\dim (W_p\cap V_q)\ge r_w(p,q),\;1\leq p,q\leq n\},$$
où la fonction de rang $r_w$ est définie par l'identité 
$r_w(p,q)=\# \{i\leq p,\; w(i)\leq q\}.$ La dimension de cette 
variété est égale à la longueur $l(w)$ de la permutation
correspondante. 

Un théorème de Lakshmibai et Sandhya permet de vérifier
combinatoirement si une variété de Schubert est ou non singulière. 

\medskip\noindent {\bf Théorème \cite{ls}.} {\em 
La variété de Schubert $X_w$ est lisse si et seulement si il 
n'existe pas de quadruplet d'entiers $i<j<k<l$ tel que }
$$w(l)<w(j)<w(k)<w(i)\qquad ou \qquad w(k)<w(l)<w(i)<w(j).$$

\medskip
Dans toute la suite, on appellera de tels quadruplets des
configurations de type 4231 et 3412, respectivement.

On sait également exprimer de façon combinatoire la dimension de
l'espace tangent de Zariski d'une variété de Schubert en chacun des 
points fixes d'un tore compatible avec le drapeau de référence. 
Ces points fixes correspondent aux permutations $v$ qui sont majorées 
par $w$ pour l'ordre de Bruhat, et la dimension de l'espace tangent de
Zariski en ce point est, d'après le résultat principal de \cite{lse}, 
égale à l'entier 
$$m(w,v) := \#\{i<j,\;\; vt_{ij}\le w\}.$$ 
Il s'ensuit que $X_v$ est incluse dans le lieu singulier de $X_w$ si
et seulement si $m(w,v)>l(w)$. 
Ceci permet en principe de déterminer ce lieu singulier, modulo
des difficultés combinatoires non négligeables. L'objet de cette note
est de décrire directement les composantes irréductibles du lieu 
singulier d'une variété de Schubert. 

Comme pouvait le laisser imaginer le théorème de Lakshmibai et
Sandhya, ces composantes sont de deux origines et correspondent à deux
types de configurations dans le diagramme de la permutation $w$. 
(Rappelons que ce diagramme est l'ensemble des points $(i,w(i))$.
Nous représenterons graphiquement ce diagramme en numérotant 
les lignes de haut en bas, et les colonnes de gauche à droite. 
Voir par exemple \cite{man} pour ces représentations.) Une
configuration 4231 sera dite minimale si elle possède les propriétés
résumées dans le diagramme ci-dessous : les zones où figurent un
symbole $\emptyset$ ne doivent pas contenir de point du diagramme; 
dans les zones où apparaît une flèche SO-NE (Sud-Ouest - Nord-Est, une
fois pour toutes), les points du
diagramme doivent être disposés selon cette direction.
A chaque choix d'un point sur chacune des deux flèches 
correspond une configuration 4231, mais on ne distingue pas entre
toutes ces configurations. 

$$
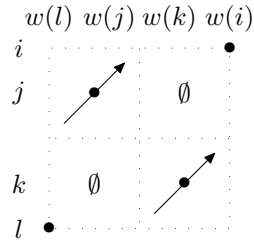
\epsfbox{sch.1}$$

De même, une configuration 3421 sera dite minimale si elle possède 
les propriétés résumées dans l'un des diagrammes ci-dessous. 
Au contraire du cas précédent, il est ici possible que n'apparaisse
aucun point du diagramme dans l'une ou l'autre des zones occupées par
une flèche.  
\smallskip

$$
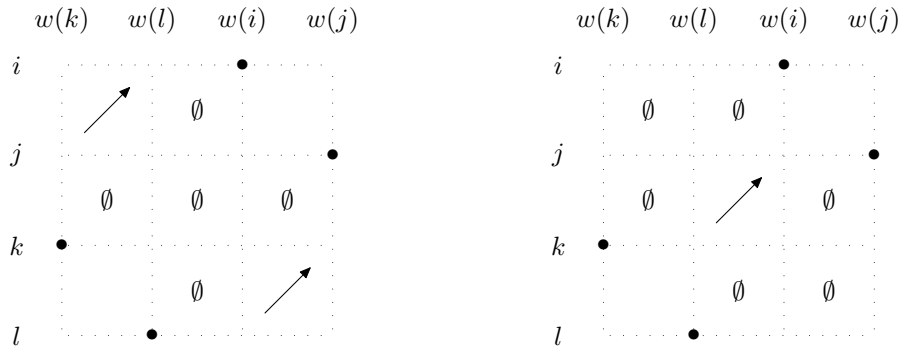
\epsfbox{sch.2}$$

\smallskip
A ces deux types de configurations minimales, on associe des
permutations  en modifiant localement le graphe
de $w$ de la façon suivante. Pour les configurations 4231 : \smallskip 

$$
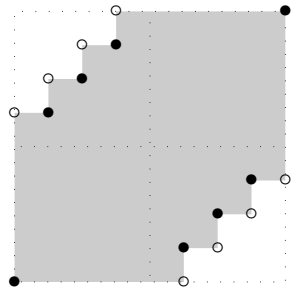
\epsfbox{sch.3}$$

\smallskip 
Pour les configurations 3412, respectivement :
\smallskip 

$$
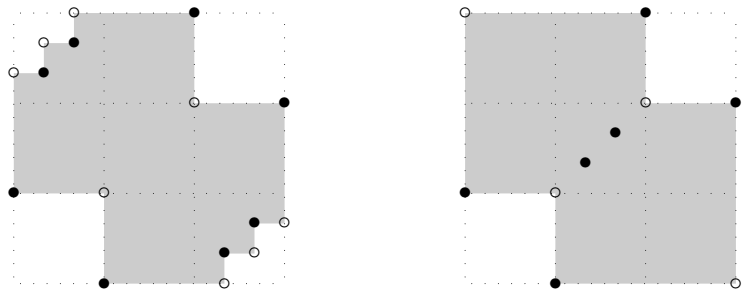
\epsfbox{sch.4}$$

\smallskip 
Dans ces figures, les $\bullet$ désignent des points du diagramme
de $w$ : on les remplace par les points désignés par des $\circ$ pour
obtenir le diagrame d'une permutation $v$. On notera que le passage 
de $w$ à $v$ augmente la fonction de rang d'une unité exactement sur
la zone grisée des figures ci-dessus, et ne la modifie pas ailleurs.
(Pour être tout à fait correct, la zone sur laquelle la fonction 
de rang augmente 
est la région polygonale figurée comme ci-dessus, mais dont
les sommets SE sont des coins externes, c'est-à-dire ne sont pas 
contenus cdans la région en question. On devra garder cette précision
à l'esprit, mais on n'en tiendra pas compte dans les figures qui suivront.)
Cela sera important pour la suite, puisque l'on utilisera
systématiquement la fonction de rang comme substitut de l'ordre de
Bruhat : on sait en effet que $v\le w$ équivaut à $r_v\ge r_w$ 
(\cite{man}, Proposition 2.1.12).

Ce sont précisément les permutations ainsi obtenues qui décrivent 
le lieu singulier des variétés de Schubert :

\begin{theo}
Les composantes irréductibles de la \vs  $X_w$ sont les \vss $X_v$,
où $v$ décrit l'ensemble $C(w)$ des permutations déduites par le procédé
ci-dessus des configurations minimales de type 4231 ou 3412 de $w$. 
\end{theo}

Le cas des variétés de Schubert des variétés de drapeaux incomplets
se déduit sans mal de cet énoncé : elles sont en effet indexées
par un ensemble restreint de permutations, de façon telle que les 
\vss correspondantes dans $\FF_n$ en soient des fibrations localement 
triviales. On peut ainsi retrouver l'énoncé classique de Svanes 
(\cite{sv}, voir aussi \cite{man}, Théorème 3.4.4) 
pour les singularités des \vss des grassmanniennes, dont les 
composantes irréductibles proviennent d'ailleurs exclusivement de
configurations de type 4231. 

\smallskip Le plan de cet article est le suivant. Dans la section 2, 
on relie les variations de l'entier $m(w,u)$ lorsque $u$ varie, à
l'existence de certaines configurations (Proposition 2). On en déduit 
dans la section 3 que les \vs $X_u$, pour $u\in C(w)$, sont incluses 
dans le lieu singulier de  $X_w$ , et maximales pour cette propriété
(Propositions 3 et 4). Dans la section 4, on montre que toute \vs 
$X_v$ incluse dans le lieu singulier de $X_w$ doit être incluse dans 
une \vs $X_u$, avec $u\in C(w)$ (Proposition 5). Le théorème 1 s'ensuit
immédiatement. Dans la section 5, enfin, on démontre une conjecture
de Lakshmibai et Sandhya, qui donne une caractérisation indirecte de
l'ensemble $C(w)$ (Théorème 8). Un cas particulier de cette conjecture
avait été obtenu relativement récemment 
par Gonciulea et Lakshmibai \cite{gl}.

\section{Préliminaires}

On a vu que l'entier $m(w,v) := \#\{i<j,\;\; vt_{ij}\le w\}$ n'était
autre que la dimension de l'espace tangent de Zariski de $X_w$ en un 
point de la cellule de Schubert associée à $v$ -- cellule qui est  
contenue dans $X_w$ si et seulement si $v\le w$ pour l'ordre de
Bruhat. La dimension de l'espace tangent de Zariski étant une fonction
semi-continue supérieurement, on doit avoir $m(w,u)\ge m(w,v)$ si
$u\le v\le w$, inégalité qui, combinatoirement parlant, ne semble pas 
évidente a priori. Comme les composantes maximales du lieu singulier
de $X_w$ seront repérées par les premières variations strictement 
positives de l'entier $m(w,v)$, il nous sera essentiel de comprendre
pourquoi et comment il varie. On notera
$$S(w,v)=\{i<j,\;\; vt_{ij}\le w\}.$$

Supposons $u=vt_{kl}$ avec $k<l$ et $v(k)>v(l)$.

On dira que $i<j$ et $k<l$ forment une configuration de type A si $v(i)<v(j)$,
et si les conditions suivantes sont vérifiées. Appelons $R_1$ le
rectangle de sommets $(i,v(i))$ (interne) et $(j,v(j))$ (externe),
appelons $R_2$ celui de sommets (externes) $(k,v(k))$ et $(l,v(l))$. 
On demande alors que 
\begin{enumerate}
\item $r_v>r_w$ sur $R_1\backslash R_2$, 
\item $r_v-r_w$ s'annule au moins une fois sur $R_1\cap R_2$. 
\end{enumerate}
Une telle configuration est représentée dans le diagramme ci-dessous
lorsque $i<k<j<l$ :
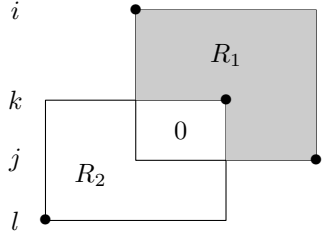
\smallskip

$$\epsfbox{sch.5}$$. 

On dira par ailleurs que $i<k<l$ (ou $k<l<i$)
forment une configuration de type $B$ 
si l'on est dans la situation suivante (ou la situation symétrique 
par rapport à la diagonale SO-NE) : $v(i)<v(l)<v(k)$, et les conditions
suivantes sont vérifiées. Appelons $R_1$ le
rectangle de sommets $(i,v(i))$ (interne) et $(k,v(l))$ (externe),
appelons $R_2$ celui de sommets $(i,v(l))$ (interne) et $(k,v(k))$
(externe), appelons enfin $R'_2$ celui de sommets $(k,v(i))$ (interne) 
et $(l,v(l))$ (externe). On demande alors que 
\begin{enumerate}
\item $r_v>r_w$ sur $R_1$, 
\item $r_v-r_w$ s'annule au moins une fois sur $R_2$ et sur $R'_2$. 
\end{enumerate}
Une telle configuration est représentée dans le diagramme ci-dessous :
\smallskip

$$
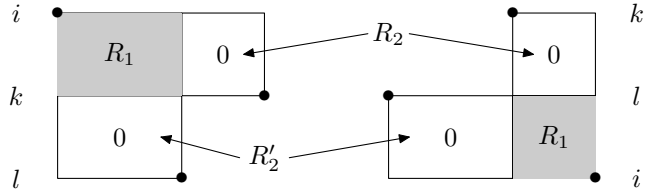
\epsfbox{sch.6}$$. 

\begin{prop} Avec les notations précédentes, la différence 
$m(w,u)-m(w,v)$ est égale
au nombre de configurations de type A ou B. \end{prop}

\dem Soient $i<j$ des entiers. On cherche à comparer les positions
de $ut_{ij}$ et $vt_{ij}$ par rapport à $w$. 

\smallskip\noindent {\em Premier cas.} Supposons que $\{i,j\}
\cap \{k,l\}=\emptyset$. Notons $R_1$ le rectangle de sommets 
$(i,v(i))=(i,u(i))$ et $(j,v(j))=(j,u(j))$, et $R_2$ le rectangle
de sommets $(k,v(k))=(k,u(l))$ et $(l,v(l))=(l,u(k))$. 
Alors les fonctions de rang de $u, v, ut_{ij}$ et $vt_{ij}$se
comparent de la façon suivante : 
$$\begin{array}{rcl}
r_u & = & r_v+\chi_{R_2}, \\
r_{ut_{ij}} & = & r_u+\e (v(i)-v(j))\chi_{R_1}, \\
r_{vt_{ij}} & = & r_v+\e (v(i)-v(j))\chi_{R_1}, 
\end{array}$$
où $\e(t)$ désigne le signe de l'entier  non nul $t$, et $\chi_R$
la fonction caractéristique du rectangle $R$. En particulier, le 
fait que $u\le v$ assure que $ut_{ij} \le vt_{ij}$. En conséquence, 
si les contributions du couple $i,j$ à $S(w,u)$ et $S(w,v)$ sont 
différentes, on doit avoir $ut_{ij}\le w$ mais pas $vt_{ij}\le w$. 
Ceci impose que $v(i)<v(j)$, et $r_u=r_v+\chi_{R_2}-\chi_{R_1}\ge
r_w$, donc $r_v>r_w$ sur $R_1\backslash R_2$. Par contre, 
$r_{vt_{ij}}=r_v-\chi_{R_1}\ngeq r_w$, ce qui signifie que $r_v-r_w$
doit s'annuler quelque part sur $R_2$, et c'est nécessairement sur 
$R_1\cap R_2$ d'après ce qui précède. On obtient donc une contribution 
à $m(w,u)-m(w,v)$ égale au nombre de configurations de type A. 

\smallskip\noindent {\em Deuxième cas.} Supposons que $\{i,j\}
\cap \{k,l\}\neq\emptyset$. Dans le cas où ces ensembles coîncident,
on a $ut_{ij})=v$ et $vt_{ij})=u$, qui sont toutes les deux majorées
par $w$. On peut donc supposer que $\{i,j\}\cap \{k,l\}$ est un 
singleton : il y a alors six cas possibles pour la position
relative de ces quatre entiers. Pour chacun de ces cas, il faut
en distinguer trois autres selon la position de $v(i)$ (ou $v(j)$)
par rapport à $v(k)$ et $v(l)$. Examinons soigneusement chacun de
ces dix-huit cas. Ils sont discutés dans les figures ci-dessous,
où l'on trouvera : à gauche, la position relative des entiers 
$i, j, k, l$; puis trois doublets de figures correspondant aux trois
positions possibles de $v(i)$ (ou $v(j)$) par rapport à $v(k)$ et
$v(l)$; chaque doublet donne à gauche les positions des points
pertinents du diagramme de $v$ (représentés par des $\bullet$)
et $vt_{ij}$ (représentés par des $\circ$, pour ceux qui sont
différents des précédents), de même à droite celles
de $u$ et $ut_{ij}$; sur chaque figure, est de plus représentée 
l'écart de la fonction de rang de $vt_{ij}$ à gauche, et $ut_{ij}$
à droite, par rapport à la fonction de rang de $v$ : à une zone gris
clair correspond un écart de $-1$, à une zone gris foncé un écart 
de $+1$. 

$$
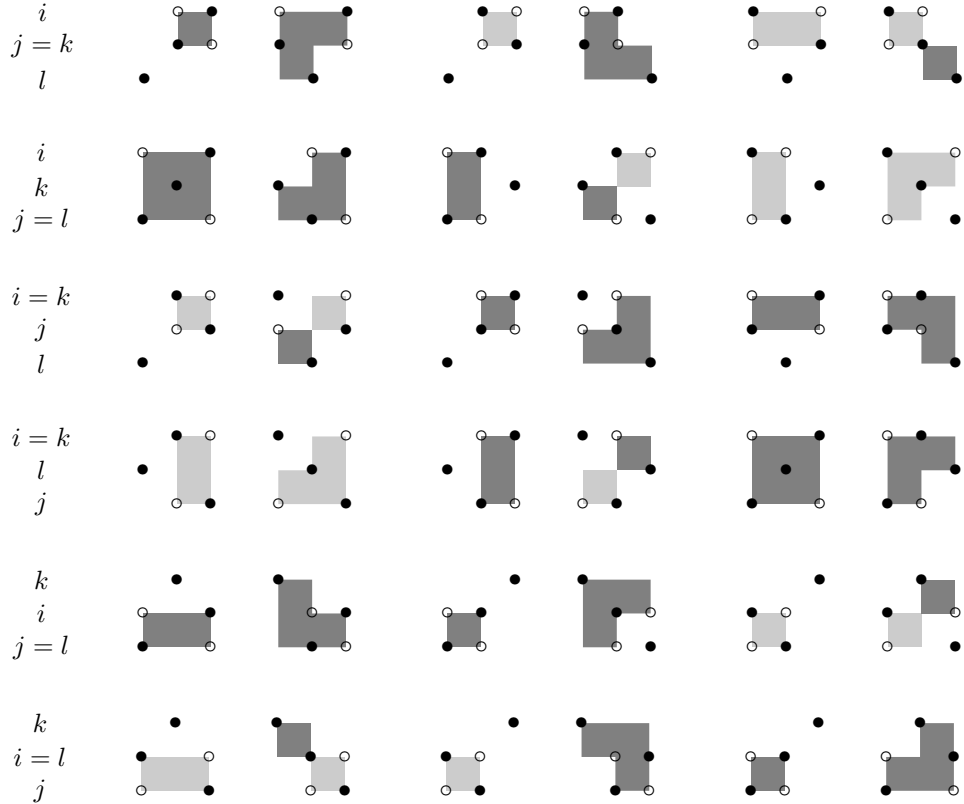
\epsfbox{sch.10}$$

$$\epsfbox{sch.11}$$

$$\epsfbox{sch.12}$$

$$\epsfbox{sch.13}$$

$$\epsfbox{sch.14}$$

$$\epsfbox{sch.15}$$

Repérons les différents cas ci-dessus par leur numéro de ligne (de 1 à
6 de haut en bas) et leur position dans la ligne (a,b ou c de gauche à
droite). On constate que pour 1a, 2a, 3b, 3c, 4c, 5a, 5b et 6c, les 
permutations $vt_{ij}$ et $ut_{ij}$ sont simultanément majorées par 
$w$. Dans les cas 3a et 5c, on a $ut_{ij}\le vt_{ij}$, et si 
$ut_{ij}\le w$, alors nécessairement $vt_{ij}\le w$. 

Les cas où $vt_{ij}\le w$ mais $ut_{ij}\nleq w$ proviennent des cas
2b, 2c, 4a et 4b, et sont symbolisés par les diagrammes ci-dessous, 
où $r_v-r_w$ doit être strictement positive sur les zones grisées, 
et doit s'annuler dans un des rectangles où figure un zéro; les
positions des points du diagramme de $v$ sur les lignes $k$ et $l$ 
sont figurées par des $\bullet$, l'autre point pertinent par un
$\circ$. 

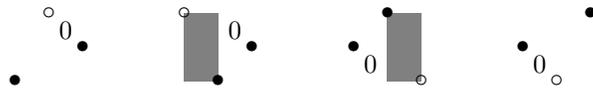
$$\epsfbox{sch.16}$$

De même, les cas où $ut_{ij}\le w$ mais $vt_{ij}\nleq w$ proviennent 
des cas 1b, 1c, 6a et 6b, et sont symbolisés par les diagrammes ci-dessous :

$$
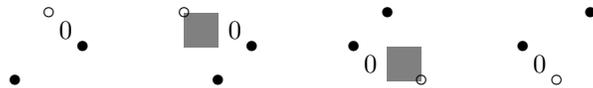
\epsfbox{sch.17}$$

On constate que les contributions de 1b et 2b d'une part, 4b et 6b 
d'autre part, se compensent exactement. Par contre, celle de 1c est 
nécessairement supérieure à celle de 2c, de même que celle de 6a est
supérieure à celle de 4a. Les différences de ces contributions
correspondent aux diagrammes suivants :

$$
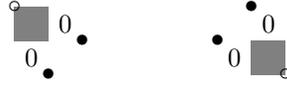
\epsfbox{sch.18}$$

On trouve exactement les configurations de type B. \qed

\section{Les composantes irréductibles du lieu singulier}

\begin{prop} Les \vss $X_u$, où $u\in C(w)$ provient d'une
  configuration minimale de type $4231$, sont des composantes
irréductibles du lieu singulier de $X_w$. \end{prop} 

\dem Montrons d'abord que $X_u$ est incluse dans le lieu singulier 
de $X_w$ : il s'agit de vérifier que $m(w,u)>l(w)$. Notons pour cela
que $u$ se déduit de $w$ par la succession de 2 séries de mouvements
élémentaires dont chacune correspond à un cycle (on appelle ici
mouvement élémentaire l'échange de deux points du diagramme qui fait
diminuer la longueur d'exactement une unité) : on commence par 
abaisser le point du diagramme en position 4 par des mouvements
élémentaires impliquant les points de la suite NO, puis on le déplace  
vers la gauche par des mouvements élémentaires impliquant les points 
de la suite SE. Les résultats de ces deux  séries de mouvements 
élémentaires sont indiqués dans la figure ci-dessous. \smallskip

$$
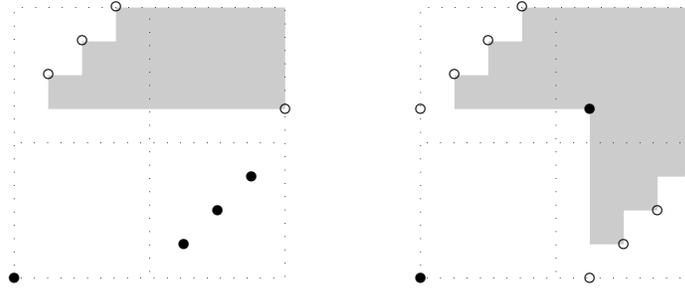
\epsfbox{sch.7}$$

\smallskip On obtient finalement une
permutation $v$, dont $u$ se déduit par un mouvement élémentaire 
impliquant le point du diagramme en position 1 et celui où a été amené 
le point initialement en position 4. Mais apparaissent ici des
configurations de type A : il suffit de prendre un point de la suite $NO$
et un point de la suite $SE$, le rectangle $R_1$ dont ils sont les sommets
satisfait les conditions qui définissent une configuration de type A.
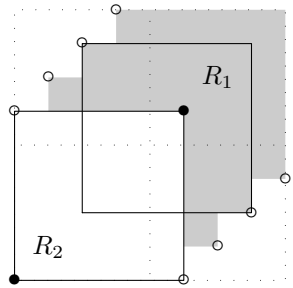
\smallskip

$$\epsfbox{sch.8}$$

\smallskip
On constate même que ce sont les seules configurations de ce type, et
qu'il n'y a pas de configuration de type B, ce dont on déduit que 
$$m(w,u)=m(w,v)+(l+1)(m+1)>m(w,v)\ge l(w).$$

Montrons maintenant que $X_u$ est une composante irréductible du 
lieu singulier de $X_w$, ce pourquoi il s'agit de démontrer que si 
$u\le x\le w$ et $l(x)=l(u)+1$, alors $m(w,x)=l(w)$. Parmi ces 
permutations $x$ figure celle que nous venons d'appeler $v$. Or on 
constate sans difficulté que $X_v$ n'est pas incluse dans le lieu
singulier de $X_w$ : en effet, dans la suite de mouvements ci-dessus
permettant de passer de $w$ à $v$, on vérifie à chaque étape que
n'apparaît aucune configuration de type A ou B. Cela implique 
au passage que $l(w)=m(w,u)-(l+1)(m+1)$, et l'on doit donc vérifier
que pour toute autre permutation $x$ comme ci-dessus, on a 
$m(w,u)=m(w,x)+(l+1)(m+1)$.

Or une telle permutation $x$ s'obtient en choisissant un point de la
suite SE, un point de la suite NO, et en les permutant, ce qui
donne la figure suivante. 

$$
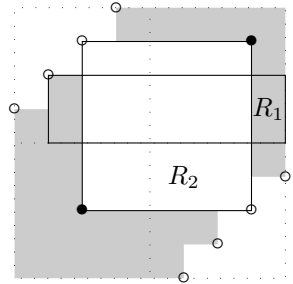
\epsfbox{sch.9}$$

On constate qu'apparaît ici une configuration de type A pour chaque 
choix d'un autre point de la suite SE, ou du point $1'$, 
et d'un autre point de la suite NO, ou du point $2'$. Donc 
$m(w,u)\ge m(w,x)+(l+1)(m+1)$, et l'égalité découle du fait 
que $m(w,x)\ge l(w)$ (ou de la constatation qu'il n'apparaît pas
d'autre configuration de type A ou B que les précédentes). 
La proposition est démontrée. \qed

\begin{prop} Les \vss $X_u$, où $u\in C(w)$ provient d'une
  configuration minimale de type $3412$, sont des composantes
irréductibles du lieu singulier de $X_w$. \end{prop} 

\dem On suit le même schéma de démonstration que pour la proposition
précédente. Rappelons que nous avons identifié deux variétés de
composantes minimales de type 3412, ce qui nous oblige à examiner
deux cas. 

\noindent {\em Premier cas.} On passe de $w$ à $u$ de la façon suivante :
on commence par abaisser le point du diagramme en position 3 par 
des mouvements élémentaires impliquant les points de la suite NO,
y compris le point du diagramme en position 1; on le déplace 
ensuite à droite à l'aide du point occupant la position 4; 
enfin, on le déplace à nouveau vers la gauche par des mouvements
élémentaires impliquant les points de la suite SE.
%, y compris le point occupant la position 2. 
Les résultats de ces trois étapes 
sont indiqués dans la figure ci-dessous. \smallskip

$$
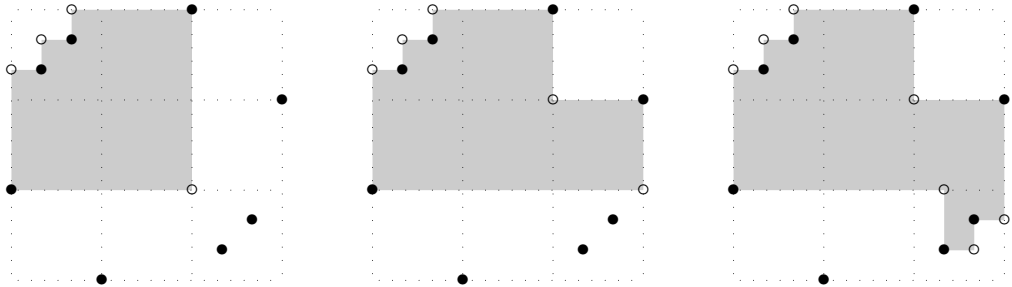
\epsfbox{schu.1}$$

\smallskip On obtient finalement une
permutation $v$, dont $u$ se déduit par un mouvement élémentaire 
impliquant le point du diagramme en position 2 et celui où a été amené 
le point initialement en position 3. Mais apparaissent ici des
configurations de type A et de type B. Pour les premières, il suffit 
de prendre un point de la suite SE : avec le point en position 33
(i.e. sur la troisième ligne en partant du haut, et la troisième
colonne en partant de la gauche, sur le quadrillage en pointillé dans
la figure), 
on obtient deux sommets opposés d'un rectangle $R_1$ satisfaisant les 
conditions qui définissent une configuration de type A. Pour les 
secondes, il suffit de prendre un point de la suite NO, y compris le
point se trouvant sur l'abscisse du point occupant initialement la 
position 1.
\smallskip

$$
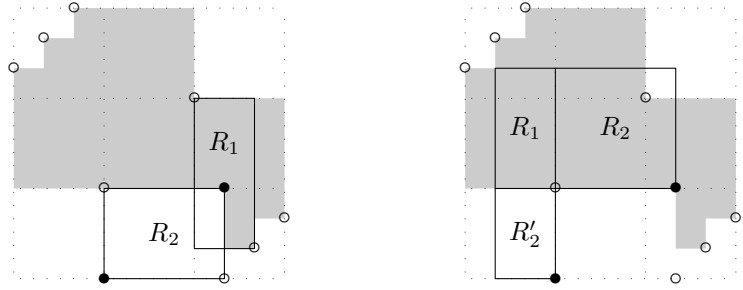
\epsfbox{schu.2}$$

\smallskip
On constate qu'il n'y a pas d'autre configuration de
type A ou B, ce qui implique que 
$$m(w,u)=m(w,v)+l+m+1>m(w,v)\ge l(w).$$
Comme dans la démonstration de la proposition précédente, on s'assure 
que l'on a en fait $m(w,v)=l(w).$

Montrons maintenant que $X_u$ est une composante irréductible du 
lieu singulier de $X_w$, ce pourquoi il s'agit de démontrer que si 
$u\le x\le w$ et $l(x)=l(u)+1$, alors $m(w,u)=m(w,x)+l+m+1$.
Or une telle permutation $x$, à symétrie près par rapport à l'une des
deux diagonales, s'obtient en choisissant un point de la SE
et en le permutant avec le point occupant la position 22. 

$$
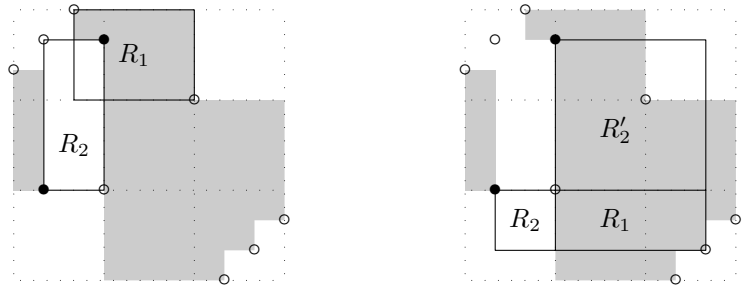
\epsfbox{schu.3}$$

On constate qu'apparaît ici une configuration de type A pour chaque 
choix d'un autre point de la suite SE (ou du point $1'$), qui est 
avec le point occupant la position 33 sommet d'un rectangle $R_1$ 
convenable. De même, apparaît une configuration de type B pour chaque
point de la suite SE. On obtient ainsi $l$ configurations de type A
et $m+1$ de type B. Donc 
$m(w,u)\ge m(w,x)+l+m+1$, et l'égalité découle du fait 
que $m(w,x)\ge l(w)$ (ou de la constatation qu'il n'apparaît pas
d'autre configuration de type A ou B que les précédentes). 

\medskip\noindent {\em Deuxième cas.} 
Ici, on passe de $w$ à $u$ de la façon suivante :
on déplace d'abord vers la gauche le point du diagramme en position 3 par 
des mouvements élémentaires impliquant les points de la suite centrale,
y compris le point du diagramme en position 1; on déplace 
ensuite le point du diagramme en position 2 vers la droite 
par des mouvements élémentaires impliquant les points de la suite
centrale (décalée).
Les résultats de ces deux étapes 
sont indiqués dans la figure ci-dessous. \smallskip

$$
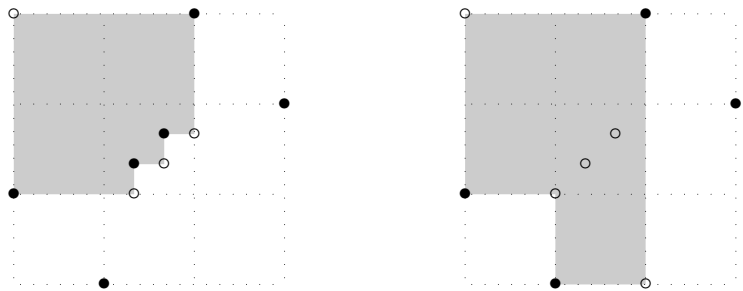
\epsfbox{schu.4}$$

\smallskip On obtient finalement une
permutation $v$, dont $u$ se déduit par un mouvement élémentaire 
impliquant le point du diagramme en position 4 et celui où a été amené 
le point initialement en position 2. Apparait  ici une et une seule
configuration de type B, correspondant au point qui occupe la position
14. \smallskip

$$
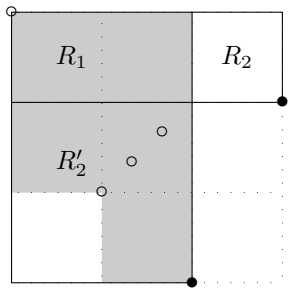
\epsfbox{schu.5}$$

\smallskip
Par contre, il n'existe aucune configuration de type A. On obtient donc
$$m(w,u)=m(w,v)+1>m(w,v)=l(w).$$

Montrons maintenant que $X_u$ est une composante irréductible du 
lieu singulier de $X_w$, ce pourquoi il s'agit de démontrer que si 
$u\le x\le w$ et $l(x)=l(u)+1$, alors $m(w,u)=m(w,x)+1$.
Or une telle permutation $x$, à symétrie près par rapport à 
la diagonale SO-NE, s'obtient en choisissant un point de la suite
centrale et en le permutant avec le point occupant la position 41. 

$$
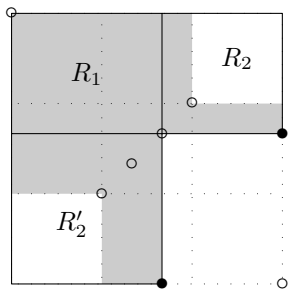
\epsfbox{schu.6}$$

On constate qu'apparaît alors une configuration de type B, 
correspondant là encore au point occupant la position 41.
La proposition est donc démontrée. \qed
 
\section{Maximalité}

Les deux propositions précédentes permettent d'affirmer que les 
variétés de Schubert $X_u$, lorsque $u$ décrit $C(w)$, sont bien des
composantes irréductibles du lieu singulier de $X_w$. Reste à vérifier
qu'il n'y en a pas d'autre, ce qui est le point le plus délicat de la
démonstration. 

\begin{prop} Si la \vs $X_x$ est incluse dans le lieu singulier de
$X_w$, alors il existe une permutation $u\in C(w)$ telle que $x\le u$.
\end{prop}

\dem On peut supposer qu'il existe une permutation $y$ telle que 
$l(y)=l(x)+1$, $x\le y\le w$, et $m(w,x)>m(w,y)$. D'après la
proposition 2, le graphe de $x$ fait alors apparaître une
configuration de type A ou B. Notons $D(x)$ l'ensemble des points
où la différence de fonctions de rang $r_x-r_w$ est strictement
positive. Si l'on choisit un chemin de $w$ à $x$, décroissant pour 
l'ordre de Bruhat et faisant diminuer la longueur d'un en un, on voit
apparaître cette région $D(x)$ comme la réunion d'une famille $\cR$
de rectangles dont les sommets sont occupés par des points du
diagramme de certaines des permutations du chemin choisi. 

Une configuration de type A ou B met en évidence des zones
contenues dans $D(x)$ (ou $D(y)$), et des points qui n'en font pas
partie. Notre stratégie consistera à isoler, à partir de ces zones
et des rectangles de la famille $\cR$ qui doivent les recouvrir, 
des configurations particulières de type 4231 ou 3412. Nous les
appellerons configurations gagnantes virtuelles, dans la mesure où 
elles ne seront a priori pas constituées de points du diagramme d'une
meme permutation intermédiaire entre $x$ et $w$. Mais nous vérifierons
que ces configurations persistent lorsque l'on ``remonte'' de $x$ à
$w$, donnant ainsi finalement une configuration du type attendu dans 
le diagramme de $w$. La position de la configuration obtenue par
rapport à $D(x)$ permettra de conclure.

\medskip\noindent {\em Premier cas.}  Partons tout d'abord d'une
configuration de type B. Dans le rectangle $R_2$, on choisit 
un point $P$ (resp. $Q$) où $r_x=r_w$, le plus à gauche (resp. 
le plus bas) possible. On choisit de façon symétrique des points 
$P'$ et $Q'$ dans $R_2'$, comme dans la figure ci-dessous. 

$$
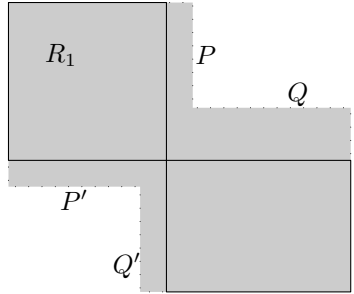
\epsfbox{schw.1}$$

On construit alors une configuration virtuelle de la façon suivante. 
Le point immédiatemment à gauche de $P$ est dans $D(x)$ par hypothèse,
il appartient donc à l'un des rectangles $R_1$ de $\cR$, 
rectangle dont le sommet NE a même abscisse que $P$. Considérons 
ensuite le point $P_1$ situé à même hauteur que $P$, immédiatement à
gauche de $R_1$. Si son abscisse est supérieure à celle de $P'$, 
soit $R_2$ un rectangle de $\cR$ contenant $P_1$. En procédant 
de cette façon, on construit une suite de sommets NE de rectangles
de $\cR$ comme dans la figure ci-dessous : à chaque étape on doit 
distinguer deux cas, selon que le sommet NE du rectangle $R_i$ 
se trouve à une hauteur inférieure au dernier point de la suite 
construite, ou non; dans le premier cas, c'est un nouveau point de
la suite, alors que dans le second cas, on ne le prend pas en compte. 

$$
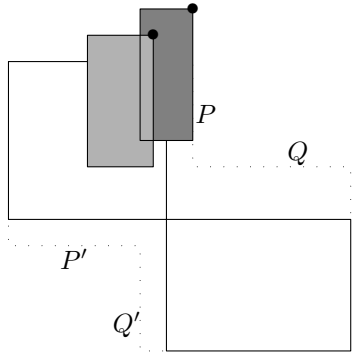
\epsfbox{schw.2}$$

On finit par obtenir un rectangle $R_m$ se trouvant au-dessus 
du point $P'$. Considérons alors le point $P'_1$ ayant même abscisse
que $P'$ et se trouvant immédiatement sous $R_m$. Si ce n'est pas
$P'$, on choisit un rectangle $R'_1$ de $\cR$ le contenant. Comme
précedemment, on construit ainsi une suite de sommets SO de rectangles
de $\cR$ situés au-dessus de $P'$. 

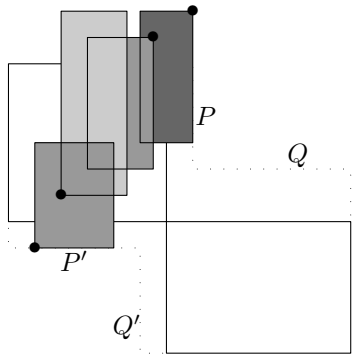
$$\epsfbox{schw.3}$$

En procédant symétriquement avec les points $Q$ et $Q'$, on obtient 
une configuration virtuelle comme ci-dessous, formée de sommets de
rectangles de $\cR$, la zone grisée étant toute entière incluse
dans $D(x)$. Apparaît en particulier une configuration 3412. 

Un point essentiel, dans cette configuration, est que chaque point 
de la suite joignant les points en position 1 et 3 est situé soit à
une altitude supérieure à celle de $P$, soit sur une abscisse
inférieure à celle de $P'$. De même, chaque point 
de la suite joignant les points en position 2 et 4 est situé soit à
une altitude inférieure à celle de $Q'$, soit sur une abscisse
supérieure à celle de $Q$. 

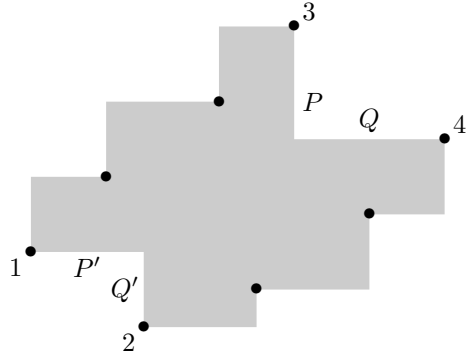
$$\epsfbox{schw.4}$$

Montrons maintenant qu'une configuration possédant
ces propriétés persiste nécessairement quand on ``remonte le temps'', 
c'est-à-dire quand on effectue une transposition qui diminue la 
longueur. Il suffit évidemment de considérer une transposition
affectant au moins l'un des points de la configuration. Si elle en
affecte deux, la présence des points $P,Q,P'$ et $Q'$ empêche que les 
points en position 1234 soient impliqués. On est donc dans la
situation de la figure ci-dessous et l'on retrouve bien une
configuration du même type en remplaçant les points $A$ et $B$
par $A'$ et $B'$. 

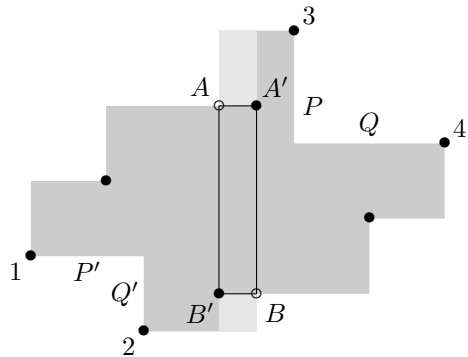
$$\epsfbox{schw.5}$$

Si la transposition n'affecte qu'un seul point de la configuration,
supposons que ce soit un des points en position 1234, qu'on peut
supposer être le point 3 du fait des symétries du problème par rapport
aux deux diagonales. A la transposition considérée correspond un
rectangle de $\cR$, donc inclus dans $D(x)$. Si le point 3 en est le
sommet SE, on le remplace par le sommet NE quitte à modifier
éventuellement la suite de points qui le joint au point 1, comme dans 
la partie gauche de la figure ci-dessous. Si le point 3 est le sommet 
NO du rectangle, qui doit alors se trouver tout entier au-dessus du
point $P$, on le remplace par le sommet SO quitte, là encore, à modifier
la suite de points qui le joint au point 1, comme dans la partie
droite de la figure ci-dessous. 

$$
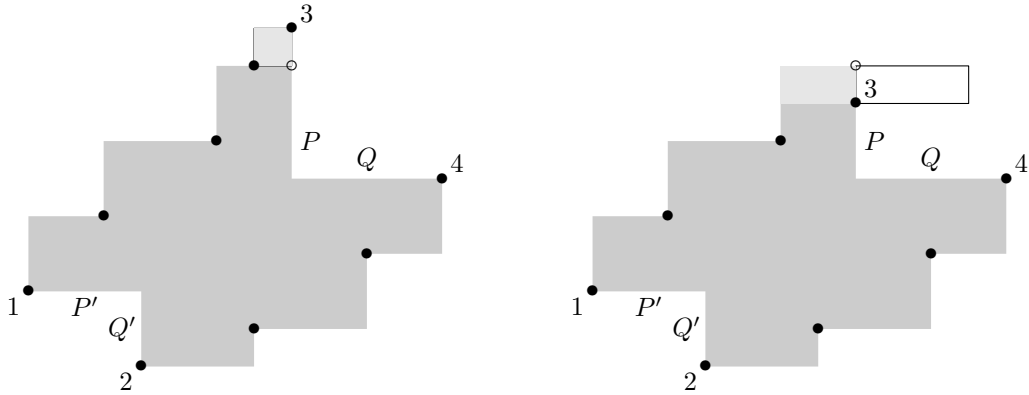
\epsfbox{schw.6}$$

Au bout du compte, on obtient une configuration possédant les mêmes 
propriétés, mais constituée de points du diagramme de $w$. On vérifie 
alors sans difficulté qu'elle doit contenir une configuration minimale
de type 3412 ou 4231. Pour cela, on peut tout d'abord supposer que la suite
joignant 1 à 3 ne possède aucun point dont l'abscisse soit incluse
entre celles de $P$ et $Q'$ : si c'est le cas, ce point doit se trouver
à une altitude supérieure à celle du point 4 et l'on peut alors le
prendre pour point 3. Au vu des symétries, on se ramène ainsi à une 
configuration du type suivant. 

$$
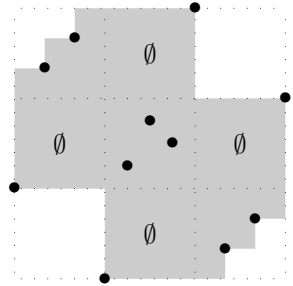
\epsfbox{schw.7}$$

Mais si l'on a des points du diagramme de $w$ à la fois dans le
rectangle central et dans le rectangle NO, par exemple, apparaît 
une configuration 4231 plus petite, et certainement minorée par une
configuration minimale de type 4231. En conséquence, soit le rectangle 
central est vide et l'on obtient une configuration minimale de type 
3412, soit il ne l'est pas et l'on peut supposer que les rectangles NO
et SE le sont. Si la suite occupée par le rectangle central est
orientée dans la direction SO-NE, on a bien affaire à une
configuration minimale de type 3412. Si ce n'est pas le cas, deux 
points orientés NO-SE de cette suite donnent, avec les points en
position 1 et 3 par exemple, une configuration 4231 plus petite, et 
certainement minorée par une configuration minimale de type 4231 
-- au sens où l'écart de sa fonction de rang avec celle de $w$ est
strictement positif sur une zone plus grande, ce qui implique bien 
qu'elle est majorée par un élément de $C(w)$.
Cela conclut l'étude de ce premier cas.

\medskip\noindent {\em Deuxième cas.} Supposons qu'il existe une
configuration de type A pour la paire $x, y$. On se placera dans
le cas où le rectangle $R_1$ se trouve au NE du rectangle $R_2$, comme
dans la figure qui précède la Proposition 2, 
les autres cas se traitant de façon
analogue (il y a deux autres cas, modulo les symétries par rapport
aux deux diagonales).  La zone $R_1 \backslash R_2$ est alors incluse
dans $D(y)$, et doit donc être couverte par une famille de rectangles 
de $\cR$,  ne contenant pas les points de $R_1\cap R_2$ qui ne sont pas
inclus dans $D(y)$ (il en existe par hypothèse). Fixons l'un de ces 
points, disons $O$. 

Supposons tout d'abord qu'il existe dans $\cR$ un rectangle, au-dessus
de $O$, et un autre à droite de $O$, dont les sommets SO soient dans 
$R_2$. Apparaît alors une configuration virtuelle 4231 comme dans la 
figure ci-dessous. 

$$
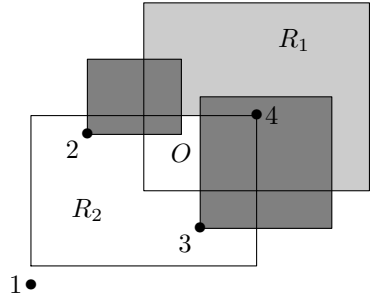
\epsfbox{schv.1}$$

Si tous les rectangles de $\cR$ se trouvant au-dessus de $O$ ont leur
sommet SO en dehors de $R_2$, on construit une configuration de la
façon suivante. On considère tout d'abord un point de $D(y)$ ayant 
même abscisse que $O$, et se trouvant au-dessus de $O$, le plus
bas possible. Ce point est contenu dans un rectangle $C_1$ de $\cR$. 
%Si ce côté supérieur de ce rectangle est au-dessous de celui de $R_2$,
On recommence alors avec le point se trouvant à même abscisse que $O$, 
immédiatement au-dessus de $C_1$. On finit par obtenir un rectangle
$C_l$ de $\cR$, dont le côté inférieur  est au-dessus du côté
supérieur de $R_2$, ou bien contenant le point $S$, sommet NO de $R_1$. 
De la suite des sommets SO de ces rectangles, plus éventuellement $S$,
on extrait alors une configuration de points comme ci-dessous.  

$$
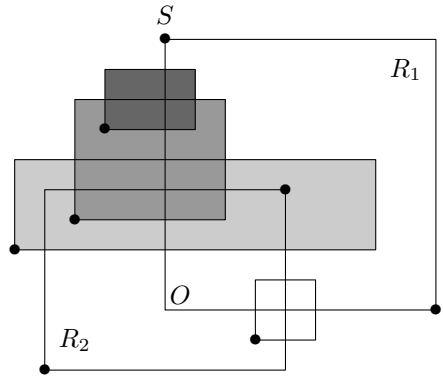
\epsfbox{schv.2}$$

Si nécessaire, on fait de même symétriquement par rapport à la
diagonale SO-NE. Dans tous les cas, on obtient une configuration du
type suivant, formé d'un ``corps'' ($C$) et de deux ``ailes''($A$ et
$A'$). Chacune 
des ailes peut  éventuellement être réduite à un seul point. 
Le point $O$ se trouve entre l'extrémité SE de l'aile NO, et 
l'extrémité NO de l'aile SE. 

$$
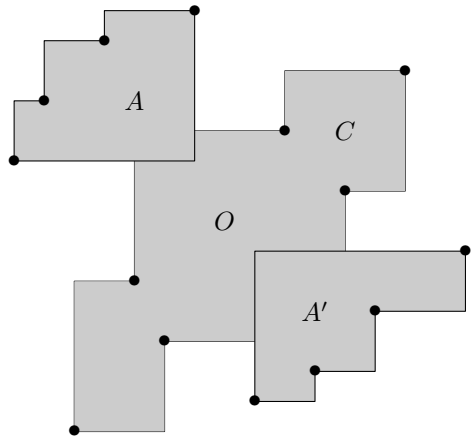
\epsfbox{schv.3}$$

L'existence du point $O$ est essentielle pour montrer que cette
configuration persiste lorsque l'on remonte de $y$ à $w$. On vérifie 
cette persistance en raisonnant au cas par cas comme dans la première
partie de cette démonstration, ce qui ne présente pas de difficulté,
le nombre de cas à examiner étant fortement limité par les symétries
de la configuration. 

Il existe donc une configuration du même type non plus virtuelle, mais
constituée de points du diagramme de $w$. Si l'une des ailes de cette
configuration n'est pas réduite à un point, on obtient une
configuration de type 3412 qui domine manifestement une configuration 
minimale du même type. Si les deux ailes sont réduites à des points, 
on obtient une configuration de type 4231 qui là encore, domine une
configuration du même type. Cela conclut la démonstration. \qed

\section{La conjecture de Lakshmibai \& Sandhya}

V. Lakshmibai et B. Sandhya conjecturent dans \cite{ls} une
caractérisation indirecte des composantes irréductibles 
du lieu singulier d'une \vs
$X_w$. Ces auteurs associent à $w$ un ensemble $Z(w)$ de permutations
$v\le w$ provenant des configurations de type 3412 ou 4231 de $w$ 
de la façon suivante. 

Si l'on part d'une configuration 3412, on demande que les points 
du diagramme de $v$ sur les colonnes correspondantes soient dans une 
configuration 1324. De plus, si l'on redresse la configuration 3412
de $w$ en une configuration 1324, la permutation obtenue $v'$ doit
vérifier $v'\le v$, alors que si l'on transforme la configuration 
1324 de $v$ en une configuration 3412, la permutation obtenue $w'$ doit
vérifier $w'\le w$. Graphiquement, cela donne la figure suivante :
à gauche, on a représenté par des $\bullet$ les points du diagramme
de $w$, par des $\circ$ ceux de $v'$, 
et à droite, par des $\bullet$ les points du diagramme
de $v$, par des $\circ$ ceux de $w'$. 
On a $r_{v'}=r_w+\chi_A$ et $r_{w'}=r_v-\chi_B$, de sorte que les 
conditions imposées à $v$ s'écrivent simplement $r_w+\chi_B\le r_v
\le r_w+\chi_A$. En particulier, la région $A$ doit contenir la 
région $B$. 

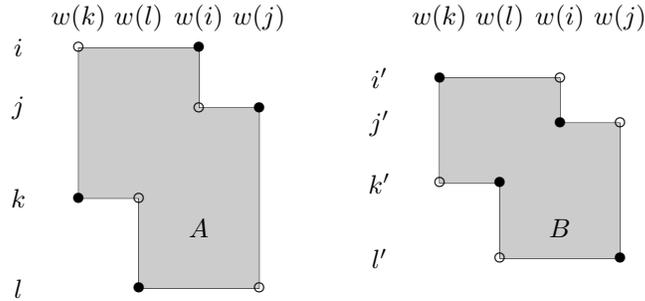
$$\epsfbox{schv.4}$$

De même, si l'on part d'une configuration 4231, on demande que les points 
du diagramme de $v$ sur les colonnes correspondantes soient dans une 
configuration 2413. De plus, si l'on redresse la configuration 4231
de $w$ en une configuration 2413, la permutation obtenue $v'$ doit
vérifier $v'\ge v$, alors que si l'on transforme la configuration 
2413 de $v$ en une configuration 4231, la permutation obtenue $w'$ doit
vérifier $w'\le w$. Graphiquement, avec les mêmes conventions que
ci-dessus, on obtient la figure suivante. Ici encore, les conditions
imposées à $v$ s'écrivent  $r_w+\chi_B\le r_v
\le r_w+\chi_A$, et la région $A$ doit contenir la région $B$. 

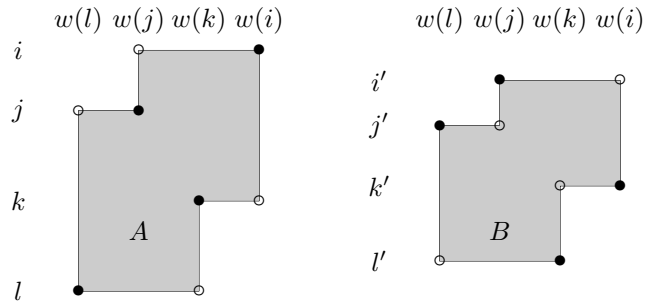
$$\epsfbox{schv.5}$$

\medskip\noindent {\bf Conjecture \cite{ls}.} {\em Les composantes
irréductibles de $X_w$ sont les \vss $X_v$, où $v$ décrit l'ensemble
des éléments maximaux de $Z(w)$ pour l'ordre de Bruhat.} 

\medskip Vérifier que le théorème 1 implique cette conjecture 
ne présente pas de difficulté sérieuse. On procède en deux temps. 

\begin{lemm} Si $u\in Z(w)$, alors le lieu singulier de $X_w$ contient 
$X_u$. \end{lemm}

\dem Supposons par exemple que $u$ provienne d'une configuration 
3412 de $w$. Comme on vient de l'expliquer, on a donc la partie 
gauche de la figure ci-dessous, où les $\bullet$ sont des points 
du diagramme de $w$ et les $\circ$ des points du diagramme de $u$. 

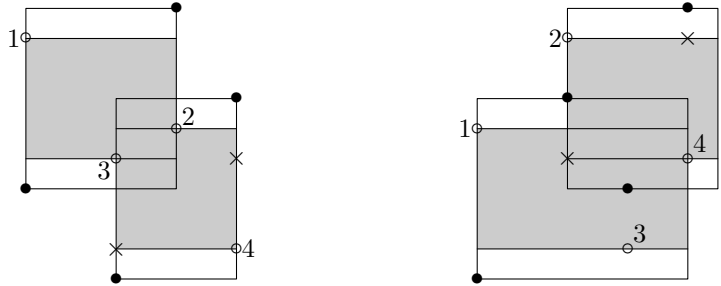
$$\epsfbox{schv.6}$$

Soit $v$ la permutation obtenue en échangeant les points du diagramme
de $u$ en position 2 et 4. On a $u\le v\le w$ : la première inégalité 
est évidente, la seconde découle de l'inégalité $r_w+\chi_B\le r_u$. 
On constate alors que le point du diagramme de $u$ (donc de $v$) en 
position 1 donne avec les deux points modifiés (indiqués par des 
$\times $ dans la figure) une configuration de 
type B. La proposition 2 implique donc que $m(w,u)>m(w,v)\ge
l(w)$, ce qui assure que $X_u$ est incluse dans le lieu singulier de 
$X_w$. 

De même, si $u$ provient d'une configuration 4231 de $w$, soit
$v$ la permutation obtenue en échangeant les points du diagramme
de $u$ en position 2 et 4, comme dans la partie droite de la 
figure ci-dessus. On a $u\le v\le w$.
On constate alors que les points du diagramme de $u$ (donc de $v$) en 
position 1 et 3 donnent avec les deux points modifiés une configuration de 
type A. La proposition 2 implique donc à nouveau que $m(w,u)>m(w,v)\ge
l(w)$, ce qui assure que $X_u$ est incluse dans le lieu singulier de 
$X_w$. \qed

\begin{lemm} L'ensemble  $C(w)$ est inclus dans $Z(w)$. \end{lemm}

\dem Considérons par exemple une permutation $v$ provenant d'une 
configuration minimale de type 4231, et choisissons une configuration 
4231 particulière de cette configuration, repérée par les indices 
$i, j, k, l$ de la figure ci-dessous. 

$$
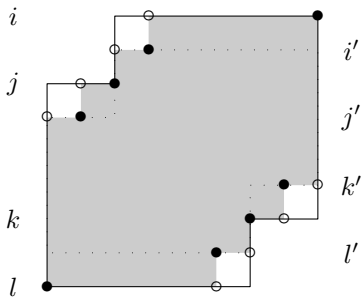
\epsfbox{schv.7}$$

On a également indiqué les indices notés ci-dessus $i', j',k',l'$. 
Les zones $A$ et $B$ s'en déduisent immédiatement : la première 
est cernée par un trait plein, la seconde par un trait pointillé dans
la figure ci-dessus. On constate
que la zone où $r_w>r_v$, qui est représentée en gris, est comprise 
entre $A$ et $B$ : c'était ce  qu'il fallait vérifier. 

On laissera le lecteur faire le raisonnement analogue pour les  
permutations provenant de configurations minimales de type  3412. \qed

\begin{theo} La conjecture de Lakshmibai et Sandhya est
  vraie. \end{theo}

\dem Il s'agit de vérifier que $C(w)$ est l'ensemble des éléments 
maximaux de $Z(w)$. Si $u$ est maximal dans $Z(w)$, le lemme 6
assure que $X_u$ est dans le lieu singulier de $X_w$, et le théorème
1  implique alors l'existence de $v\in C(w)$ tel que $u\le v$. 
Comme $v\in Z(w)$ d'après le lemme 7, et $u$ est maximal dans 
cet ensemble, il vient $u=v$. Réciproquement, si $u\in C(w)$ et $u\lneq
v$ avec $v\in Z(w)$, alors $X_v$ n'est pas incluse dans le lieu
singulier de $X_w$, ce qui contredit le lemme 6. 
\qed

\medskip\noindent
{\sc Laurent Manivel, Institut Fourier}, UMR 5582 du CNRS, Université 
Joseph Fourier, BP 74, 38402 Saint Martin d'Hères, France. 

\noindent {\em E-mail} : Laurent.Manivel@ujf-grenoble.fr

\end{document}